\documentclass[11pt]{article}
\usepackage{style}
\graphicspath{ {./images/} }

\usepackage{blindtext}
\usepackage{algorithm,algpseudocode}
\usepackage{graphicx}
\usepackage{natbib}
\usepackage{subfiles}
\usepackage{cleveref}
\usepackage{dsfont}
\usepackage{mathtools}
\usepackage{multirow}
\usepackage{authblk}
\usepackage{setspace}[]

\title{Minimax Optimal Probability Matrix Estimation For Graphon With Spectral Decay}

\begin{document}
\author[1]{Yuchen Chen}
\author[1]{Jing Lei}
\affil[1]{Carnegie Mellon University}
\maketitle

\begin{abstract}
    We study the optimal estimation of probability matrices of random graph models generated from graphons. This problem has been extensively studied in the case of step-graphons and H\"{o}lder smooth graphons. In this work, we characterize the regularity of graphons based on the decay rates of their eigenvalues. Our results show that for such classes of graphons, the minimax upper bound is achieved by a spectral thresholding algorithm and matches an information-theoretic lower bound up to a log factor. We provide insights on potential sources of this extra logarithm factor and discuss scenarios where exactly matching bounds can be obtained. This marks a difference from the step-graphon and H\"{o}lder smooth settings, because in those settings, there is a known computational-statistical gap where no polynomial time algorithm can achieve the statistical minimax rate. This contrast reflects a deeper observation that the spectral decay is an intrinsic feature of a graphon while smoothness is not.
\end{abstract}

\setstretch{1.25}

\section{Introduction}
Many modern datasets have network structures, comprising relational links among entities. For example, data on social networks may consist of friendship links between profiles on a social media service. For more examples of network problems and data, we refer the reader to \citet{girvan_community_2002,wasserman_social_1994,goldenberg_survey_2010,kolaczyk_statistical_2009,newman_networks_2010}. Network data is typically represented by an undirected graph where the nodes are the entities and the vertices are the relational links. Mathematically, we can define an undirected graph with $n$ nodes as an $n\times n$ symmetric matrix $A$, where $a_{ij} = 1$ if nodes $i$ and $j$ are connected and $a_{ij} = 0$ otherwise, called an adjacency matrix.

In statistical analysis of networks, we assume our observed graph is generated from a random graph model. A natural class of random graph models are symmetric exchangeable random graphs, which assumes the random graph has the same distribution under arbitrary relabeling or permutation of the nodes. Results from \citet{aldous_representations_1981,hoover_relations_1979} show that all infinite exchangeable random graphs can be generated by a symmetric function $W: [0,1] \times [0,1] \rightarrow [0,1]$ called a graphon. Such random graphs are formed by first sampling $\xi_1,\xi_2,\dots$ independently from Uniform$([0,1])$. Two nodes $i,j$ are connected with probability $W(\xi_i,\xi_j)$.  This result is the two-way counterpart of the celebrated de Finetti's theorem. Intuitively, the nodes in an exchangeable random graph serves as a sample point randomly drawn from a population. But instead of measuring attributes only about this individual as in traditional iid sampling, a random graph records the interaction between pairs of individuals. We note that graphons were first introduced through the study of graph limits. \citet{lovasz_large_2012} is a comprehensive reference for this viewpoint.

To get an $n$-node network, we can truncate this infinite array. That is we first sample latent node variables $\xi_1,...,\xi_n \sim_{iid} \text{Uniform}([0,1])$. Using these latent node variables, we can form the probability matrix $M$ where
\begin{equation}\label{eq:W2M}
    M_{ij} = W(\xi_i,\xi_j)\mathds{1}(i\neq j),
    \end{equation}
for $i,j=1,...,n$. Then the adjacency matrix $A$ is generated by
$$
    A_{ij} \sim \text{Bernoulli}(M_{ij})\,,~\text{independently for}~~ 1\le i<j\le n\,.
$$

In this setting, there is a natural identification issue. Two different graphons may induce the same random graph distribution. We say two graphons $W_1,W_2$ are equivalent if there exists measure-preserving mappings $h_1,h_2:[0,1]\mapsto [0,1]$ such that
$$
    W_1(h_1(x),h_1(y)) = W_2(h_2(x),h_2(y)) \text{ a.e. } 
$$
Thus, graphons can only be identified up to equivalence class given by the above equivalence relation. To get around this issue, we focus our analysis on estimating the probability matrix $M$ corresponding to the graphon $W$ under the loss
$$
    \frac{1}{n^2} \lfro{\widehat{M}-M}^2.
$$
That is, given the observed adjacency matrix $A$, how well can we estimate the probability of a connection between two nodes? Understanding such connection probabilities is a fundamental task in statistical network analysis.  More details about the identifiability of graphons are given in \Cref{sec:equivalence}.

\paragraph{Previous Work}
The probability matrix estimation problem has been extensively studied for stochastic block models and H\"{o}lder smooth graphons \citep{gao_rate-optimal_2015,klopp_oracle_2017}. These works establish the minimax rate for probability matrix estimation for both the dense and sparse cases.

In both \citet{gao_rate-optimal_2015,klopp_oracle_2017}, the algorithm that achieves the minimax rate runs in exponential time and is thus not computationally feasible. Further research has been done to analyze computationally feasible methods for probability matrix estimation. A computationally friendly spectral method, universal singular value thresholding (USVT), has been proposed as an alternative method \citep{chatterjee_matrix_2015,xu_rates_2018}. However, this method does not achieve the minimax rates established in \citet{gao_rate-optimal_2015,klopp_oracle_2017}.

More recently, there have been efforts to understand the computational-statistical gap. Using the low degree polynomial framework established by \citet{schramm_computational_2022}, \citet{luo_computational_2023} establish computational lower bounds for the probability matrix problem for both stochastic block models and H\"{o}lder smooth graphons. Their work shows that indeed polynomial time algorithms cannot achieve the minimax rate. They establish a computational lower bound and show that the USVT rate for stochastic block models matches up to logarithm terms, but still the USVT rate for H\"{o}lder smooth graphons does not.

\paragraph{Contributions}
In the previous work, we see that in the stochastic block model and H\"{o}lder smooth graphon setting, there is a statistical-computational gap. The computationally efficient spectral method does not achieve the statistical minimax rate. This raises the question: is there a class of graphons where USVT is minimax optimal? In this paper, we give an affirmative answer to this question by viewing graphons as integral operators and looking at their spectral decomposition. 

By viewing graphons as integral operators, they have a spectral decomposition
$$
    W(x,y) = \sum_{i=1}^\infty \omega_i \varphi_i(x)\varphi_i(y),
$$
where $\omega_i$ are eigenvalues and $\varphi_i$ are eigenfunctions. This decomposition is analogous to the spectral decomposition of symmetric matrices. If we label the eigenvalue/function pairs such that $|\omega_1| \geq |\omega_2| \geq \cdots$, it is known that $\omega_i \rightarrow 0$ as $i \rightarrow \infty$. A natural way to characterize the regularity of $W$ is by characterizing how fast the eigenvalues converge to $0$. In this work, we consider the case where the eigenvalues decay polynomially, i.e
$$
    |\omega_k| \lesssim k^{-\alpha},
$$
for all $k$ and a fixed $\alpha$. This type of decay condition is similar to other non-parametric estimation problems such as wavelets \citep{donoho_minimax_1998} or functional data analysis \citep{hall_methodology_2007}
. The decay of eigenvalues in the eigen-decomposition of linear operators has also been extensively studied in functional analysis \citet{lax_functional_2002,konig_eigenvalue_1986,birman_estimates_1977}.

For $\alpha > 1$, we find that spectral methods are minimax optimal up to log terms with upper and lower bounds 
$$
    \sup_{\mathcal{W}(\alpha)} \frac{1}{n^2} E\left[\lfro{\widehat{M}-M}^2\right] \lesssim n^{-\frac{2\alpha -1}{2\alpha}},
$$
$$
    \inf_{\widehat{M}} \sup_{\mathcal{W}(\alpha)} \frac{1}{n^2} E\left[\lfro{\widehat{M}-M}^2\right] \gtrsim (n\log(n))^{-\frac{2\alpha -1}{2\alpha}},
$$
where $\mathcal{W}(\alpha)$ is the set of all graphons $W$ with polynomial decay of rate $\alpha$. We give a more rigorous definition of this in Section \ref{sec: main results}. Further discussions on the logarithm factor in the lower bound is given at the end of Section \ref{sec: main results}, where the source of this extra logarithm factor is likely due to the uncertainty in the latent variables $\xi_i$ and the lack of a more refined packing number bound on the Stiefel manifold subject to an $\ell_\infty$ norm constraint. Based on these insights, we provide two scenarios where the lower bound can be improved to the exact rate of $n^{-\frac{2\alpha-1}{2\alpha}}$.

In contrast, under the H\"{o}lder smooth setting, the minimax lower bound is given by $$
    \inf_{\widehat{M}}\sup_{\mathcal{H}(\gamma)}\sup_{P_\xi} \frac{1}{n^2} E\left[\lfro{\widehat{M}-M}^2\right] \gtrsim \begin{cases}
        n^{-\frac{2\gamma}{\gamma+1}} & 0<\gamma<1\\
        \frac{\log n}{n} &\gamma > 1,
    \end{cases}
$$ where $\mathcal{H}(\gamma)$ denotes the corresponding H\"{o}lder class and $P_\xi$ denotes the sampling distribution of the latent $\xi$ \citep{gao_rate-optimal_2015}. The upper bound of the USVT method is given by $$
    \sup_{\mathcal{H}(\gamma)}\sup_{P_\xi} \frac{1}{n^2} E\left[\lfro{\widehat{M}-M}^2\right] \lesssim n^{-\frac{2\gamma}{2\gamma+1}}
$$\citep{xu_rates_2018}. We note that in the previous works, the rates are computed over an arbitrary distribution $P_\xi$ for the latent $\xi$, while we focus on the setting where $P_\xi$ are iid Uniform$([0,1])$. The uniform sampling is natural from the exchangeable array and random graph viewpoint \citep{aldous_representations_1981,hoover_relations_1979,lovasz_large_2012}. If we consider arbitrary $P_\xi$, we can achieve exactly matching rates
    $$
            \inf_{\hat{M}} \sup_{\mathcal{W
            }(\alpha)} \sup_{P_\xi}\frac{1}{n^2} E\left[\lfro{\widehat{M}-M}^2\right] \asymp n^{-\frac{2\alpha-1}{2\alpha}}\,.
    $$
We believe these results bring further insight into the optimality of spectral methods and give further support for the use of spectral methods in statistical network analysis.

\paragraph{Organization}
The paper is organized as follows. In Section \ref{sec: spectral}, we define spectral properties of graphons and define the class of graphons that we will analyze in detail. In Section \ref{sec: main results}, we state the main results. Section \ref{sec: proofs} gives the proofs of the main results. More technical details in the proofs are left to Appendices \ref{sec: main appendix}, \ref{sec: matching} and \ref{sec: auxiliary appendix}.

\paragraph{Notation} For any matrix $M$, we use $\lfro{M} := \sqrt{\sum_{i,j} M_{ij}^2}$ to denote the Frobenius norm, $\norm{M}_\infty$ to denote the largest entry of $M$ and $\opnorm{M}$ to denote the operator norm of the matrix $M$. We use $L_2([0,1])$ to denote the Hilbert space of square-integrable functions on $[0,1]$ equipped with inner product $\langle f, g \rangle = \int f g$ and $\norm{f}_2 := \sqrt{\int f^2}$. For two sequences $a(n)$ and $b(n)$, we say that $b(n) \lesssim b(n)$ if there exists some constant $C$ such that $a(n) \leq Cb(n)$. We say $a(n) \asymp b(n)$ if $a(n) \lesssim b(n)$ and $a(n) \gtrsim b(n)$. In this paper we will need to denote eigenvalues of graphons when viewed as a linear operator and eigenvalues of their corresponding probability matrix. To avoid confusion, we use $\omega_i$ to denote eigenvalues of graphons and $\lambda_i$ to denote eigenvalues of probability matrices.

\section{Spectral Properties of Graphons}
\label{sec: spectral}
In this section, we give an overview of the spectral decomposition of graphons, define the notion of spectral decay in the case of graphons, and discuss the equivalence class in the graphon space.  

\subsection{An operator perspective to graphons}
The spectral theory of graphon arises from viewing graphons as linear operators on $L_2([0,1])$. This is done by viewing a graphon $W$ as the kernel of an integral kernel operator.

\begin{definition}
    Given a graphon $W$, we define the graphon operator
    $$
        (T_W X)(v) := \int_0^1 W(u,v)X(u) du
    $$
    which is a linear operator on $L_2([0,1])$.
\end{definition}

In particular, the graphon operator is a symmetric compact operator from $L_2[0,1]$ to itself, which has an eigendecomposition \citep{lax_functional_2002}. 

\begin{definition}
    We say $\varphi:[0,1]\rightarrow \R$ is an eigenfunction of $T_W$ with eigenvalue $\omega$ if 
    $$
        (T_W \varphi)(v) = \omega \varphi(v).
    $$
\end{definition}

By the spectral theorem for symmetric compact operators on a Hilbert space \citep{lax_functional_2002}[Theorem 3, Chapter 28], there exists an orthonormal basis of $L_2([0,1])$ consisting of eigenfunctions of $T_W$. Let $\{\varphi_i\}$ denote such a basis and $\omega_i$ be the corresponding eigenvalues. Then we have a spectral decomposition of $W$ given by 
\begin{equation}
\label{eq: graphon decompositon}
    W(u,v) = \sum_{i=1}^\infty \omega_i \varphi(u)\varphi(v).
\end{equation}
 The convergence in the spectral decomposition is $L_2$-convergence.

In this paper, the graphons we will work with are trace class.
\begin{definition}
    \label{def: trace-class}
    We say a graphon $W$ is trace-class if
    $$
        \sum_{i=1}^\infty |\omega_i| < \infty,
    $$
    where $\omega_i$ are the eigenvalues of $W$.
\end{definition}

A key consequence of being trace-class is that $W$ admits a strong spectral decomposition \citep{lei_network_2021}. That is for any trace-class graphon $W$ with spectral decomposition $\sum_{i=1}^\infty \omega_i \varphi(u)\varphi(v)$, we have
\begin{equation}
\label{eq: strong spectral decomposition}
    \sum_{i=1}^\infty |\omega_i| \varphi_i^2(x) < \infty \text{ a.e.}.
\end{equation}

Strong spectral decomposition further implies that Equation \ref{eq: graphon decompositon} holds almost everywhere as well \citep{lei_network_2021}. As changing the graphon on a measure-zero set does not affect the corresponding network model, we will suppose that Equation \ref{eq: graphon decompositon} holds everywhere. For additional references on trace class operators see \citet{lax_functional_2002}.

In the remainder of this paper, we will order the eigenvalues in decreasing order in absolute value, i.e. $|\omega_1| \geq |\omega_2| \geq \dots$. As $W$ is bounded between $[0,1]$, the eigenvalues will be contained in $[-1,1]$. Furthermore, they will concentrate around $0$, i.e $\lim \omega_i \rightarrow 0.$ How fast the eigenvalues decay to $0$ is a natural measure of the regularity of $W$. The faster they decay the ``simpler'' the graphon is. We want to characterize rates of estimating the probability matrix when the eigenvalues $\omega_i$ decay sufficiently fast through the following decay condition.

\begin{definition}
\label{def: spectral decay}
    We say graphon $W$ has polynomial spectral decay with rate $\alpha$ if
    $$
        |\omega_k| \leq Ck^{-\alpha},
    $$
    for all $k$ where $C$ is a constant.
\end{definition}

Such spectral decay conditions have been studied in other areas such as wavelets \citet{donoho_minimax_1998} and functional data analysis \citep{hall_methodology_2007}. Moreover, this class of graphons contains stochastic block models and smooth graphons, which have been studied previously \citep{gao_rate-optimal_2015,xu_rates_2018,klopp_oracle_2017}. We illustrate this connection in the following two examples.

\begin{example}[Block Models]
    Consider a stochastic block model [SBM]\citep{holland1983stochastic} with $n$ nodes and $k$ communities. In this model, each of the nodes is sorted into one of the $k$ communities. The connection probability is encoded by $B \in [0,1]^{k\times k}$, a symmetric $k\times k$ matrix, where $B_{ij}$ denotes the probability of a connection between the $i$th and $j$th community.

    SBMs correspond to step graphons. We can parameterize the step graphon by a membership function $z:[0,1] \rightarrow [k]$ and the matrix $B$. Then we can write the step graphon as
    $$
        W(x,y) = B_{z(x),z(y)}.
    $$
    Let's assume that the communities have equal size. That is $\int \mathds{1}(z(x)=j) = \frac{1}{k}$ for all $j=1,...,k$. In this case, the spectral decomposition of this graphon corresponds to the spectral decomposition of the matrix $B$.

    \begin{proposition}
        Let $W$ be a step graphon parameterized by membership function $z$ and probability matrix $B$ with equal-sized communities. Then the eigenvalues of $W$ are proportional to the eigenvalues of $B$.
    \end{proposition}

    \begin{proof}
        Let $B=U\Lambda U^T$ be the eigen-decomposition of $B$ and let $u_i$ denote the columns of $U$ and $\lambda_i$ denote the eigenvalues. From this decomposition, we know that $B_{st} = \sum_{i=1}^k \lambda_i u_{is}u_{it}$. Define functions
        $$
            \varphi_i(x) := \sum_{j=1}^k u_{ij}\frac{1}{\sqrt{k}} \mathds{1}(z(x) = k).
        $$
        One can check using that the eigenvectors $u_i$ are orthonormal that $\{\varphi_i\}_{i=1}^k$ form an orthonormal set in $L_2([0,1])$. The graphon $W$ has spectral decomposition
        \begin{align*}
            W(x,y) = & \sum_{i=1}^k \lambda_i k \varphi_i(x)\varphi_i(y). 
        \end{align*}
            Thus, an equal sized SBM has polynomial decay at a rate of $\alpha$ is determined by the polynomial decay of the probability matrix.
    \end{proof}

    Similar analysis shows that extensions of the SBM such as degree corrected stochastic block models or mixed membership stochastic block models can also have polynomial spectral decay.
\end{example}

\begin{example}[Smooth Graphons]
Suppose that the graphon $W$ has $b$ continuous derivatives. Then we have the following eigenvalue decay.

\begin{proposition}
    \label{prop: smooth graphon decay}
    If the graphon $W$ has $b$ continuous derivatives, then $W$ has polynomial spectral decay with rate proportional to $b$.
\end{proposition}

This result is shown in the proof of Theorem 13, Chapter 30 in \citet{lax_functional_2002} and is a combination of the Courant's minimax formulation of eigenvalues for linear operators, a functional analysis version of the minimax formulation of eigenvalues of symmetric matrices, and approximation theory. We sketch the idea below.
\begin{proof}
Using Courant's principle Theorem 4 Chapter 28 in \citet{lax_functional_2002} we know that the $i$-th eigenvalue of $W$ satisfies
$$
    \omega_i \leq \max_{u \perp S_{i-1}, \norm{u} = 1} \int W(x,y)u(x)u(y) dxdy,
$$
where $S_{i-1}$ is any subspace of dimension $i-1$. We can take $S_{i-1}$ to be the space of all polynomials in $x,y$ with degree less than $i-1$. A result in approximation theory shows that $W$ can be approximated by such a polynomial $P_i$ with error
$$
    \norm{W-P_i}_2^2 \leq C i^{-\Tilde{b}},
$$
where $\Tilde{b}$ is proportional to $b$.
Now as $u$ is orthogonal to $S_{i-1}$ and $P_i \in S_{i-1}$, see that
$$
\int W(x,y)u(x)u(y) dxdy = \int [W(x,y) -P_i(x,y)]u(x)u(y) dxdy.
$$
Now use Cauchy-Schwartz to see that
\begin{equation*}
    \begin{split}
        \int [W(x,y) -P_i(x,y)]u(x)u(y) dxdy &\leq \norm{W-P_i}_2^2 \norm{u}_2^2\norm{u}_2^2 \\
        &\leq C i^{-\Tilde{b}}.
    \end{split}
\end{equation*}
As a result, we have $\omega_i \leq C i^{-\Tilde{b}}$, showing the desired decay rate.
\end{proof}

\end{example}

Therefore, graphons with polynomial decay include but are not limited to SBM and H\"{o}lder smooth graphons studied in previous work. The decay of eigenvalues of integral operators is a well-studied topic in functional analysis; additional results of this type may be found in \citet{lax_functional_2002,birman_estimates_1977,konig_eigenvalue_1986}.

\subsection{Graphon equivalence classes}\label{sec:equivalence}
It is well-known that graphons are identifiable only up to equivalence by measure-preserving transformation \citep{lovasz_large_2012,lei_network_2021}. Suppose that $h:[0,1] \rightarrow [0,1]$ is a measure preserving bijection such that $h(A)$ has the same Lebesgue measure as $A$ for all Borel measurable sets $A$, then for any graphon $W(\cdot,\cdot)$, the graphon $W(h(\cdot),h(\cdot))$ induces the same distribution of the observed random graph as the original graphon $W(\cdot,\cdot)$.

Given the natural identifiability limitation of graphons, any property associated with a graphon is considered intrinsic if it holds simultaneously for all members in the same equivalence class, and vice versa.

\paragraph{Spectral decay is an intrinsic graphon property.}
Suppose we have a graphon $W$ with spectral decomposition $W(x,y)=\sum_{i=1}^\infty \omega_i \varphi_i(x)\varphi_i(y)$, and a measure preserving mapping $h(\cdot)$. The transformed graphon has the form $$W(h(x),h(y))=\sum_{i=1}^\infty \omega_i \varphi_i(h(x))\varphi_i(h(y))\,.$$ We claim that this is its spectral decomposition. We can see this by observing the isomorphism between the Hilbert spaces $L_2([0,1])$ and $L_2^h([0,1])$, with the latter equipped with the inner product $\langle f,g\rangle_{L_2^h([0,1])}=\int f(h(x))g(h(x))dx$, and correspondence $f(\cdot)\mapsto (f\circ h)(\cdot)\equiv f(h(\cdot))$. For example, for any $f$, $g$ in $L_2([0,1])$, 
%$\{\varphi_i \circ h\}_{i=1}^\infty$ is an orthonormal basis holds because
\begin{align*}
   \langle f\circ h,g\circ h\rangle_{L_2^h([0,1])}=& \int f(h(x)) g(h(x)) d\mu(x) = \int f(x) g(x) dh_*\mu(x)\\
   =&\int f(x) g(x) d\mu(x)=\langle f, g\rangle_{L_2([0,1])}\,,
\end{align*}
where $\mu$ is the Lebesgue measure and $h_*$ denotes the pushforward by $h$. As $h$ is measure-preserving, we know that $h_*\mu = \mu$.
% Thus, the orthogonality of $\{\varphi_i \circ h\}_{i=1}^\infty$ follows from the orthogonality of $\{\varphi_i\}_{i=1}^\infty$.

As a result, all graphons in the same equivalence class have the same sequence of eigenvalues.  So any property about the eigenvalue sequence, including the eigen decay property, is an intrinsic graphon property.

\paragraph{Smoothness is not intrinsic for graphons.}
On the other hand, smoothness may not be preserved by a measure-preserving transformation. For example, consider the graphon $W(x,y) = y$. We can break the smoothness of this graphon using a measure preserving transformation. One example of such a transformation is
$$
h(x) =  \begin{cases} 
      \frac{1}{2} + x & x\leq \frac{1}{2} \\
      x-\frac{1}{2} & x > \frac{1}{2}.
   \end{cases}
$$
Another class of measure-preserving mappings $h$ that can make the derivatives of $f\circ h$ arbitrarily large for any smooth, non-zero $f$ is $h(x)=h_n(x)=nx-\lfloor nx\rfloor$.

\subsection{Diagonal of a graphon}
The diagonal values of a graphon are not involved in the probability matrix if we assume the random graph does not have self-loops.  But they may appear when we relate the spectral properties of the graphon to those of the probability matrix.  When the graphon is assumed to be smooth, such diagonal entries can be controlled using the off-diagonal entries according to smoothness.
In order to treat general graphons, we make a mild regularity assumption on the diagonal of the graphon.
\begin{assumption}
\label{as: L_2 diagonal}
    Assume that the partial sums of the diagonal, $\sum_{i=1}^k \omega_i \varphi_i^2(x)$, converge to $W(x,x)$ in $L_2$.
\end{assumption}

We clarify why this assumption is not too restrictive. Indeed we know that the partial sums of the diagonal converge almost everywhere (w.r.t Lebesgue measure on $[0,1]$) to the diagonal of the graphon \citep{brislawn_traceable_1991}[Theorem 3.1]. We are ruling out the exceptional cases where $L_2$ convergence does not hold. Here are some conditions where Assumption \ref{as: L_2 diagonal} hold.

\begin{proposition}
    \label{prop: as 1 holds}
    If either
    \begin{enumerate}
        \item There exists a constant $B$ such that $\sum_{i=1}^k \omega_i\varphi^2_i(x) < B$ for all $k,x$.
        \item The partial sums of the diagonal converge uniformly to the diagonal of the graphon.
    \end{enumerate}
    then Assumption \ref{as: L_2 diagonal} holds.
\end{proposition}

In either of these conditions, we are excluding the graphons where the diagonals can get too chaotic. In such cases, the error of the diagonal may be too overwhelming. 

Here are some examples of graphons where Assumption \ref{as: L_2 diagonal} is satisfied.
\begin{enumerate}
    \item Finite rank graphons. In this case, it is clear that the diagonals must converge uniformly. This class includes the popular stochastic block models.
    \item Continuous graphons: It is known that if the graphon is continuous then its eigenfunctions are continuous as well. We can split the graphon $W$ into its positive and negative components $W^+$ and $W^-$ where $W^+(x,y) = \sum_i \omega_i^+ \varphi_i^+(x)\varphi_i^+(y)$ and $W^-(x,y) = \sum_i \omega_i^- \varphi_i^-(x)\varphi_i^-(y)$. Here the superscripts denote the eigenvalue/function pairs corresponding to positive eigenvalues and negative eigenvalues respectively. Then applying Dini's theorem shows that the partial sums of $W^+$ and $W^-$ converge uniformly.
    \item Trigonometric Basis: If the eigenfunctions consist of the popular trigonometric basis of $L_2([0,1])$, then Assumption $1$ is also satisfied. In this case, the partial sums are uniformly bounded.
\end{enumerate}

For the remainder of the paper we use $\mathcal{W}(\alpha)$ to denote graphons with polynomial decay $\alpha$ satisfying Assumption \ref{as: L_2 diagonal}.

\section{Main Results}
\label{sec: main results}
In this section, we present the main results of this paper. We begin by introducing the singular value thresholding algorithm that achieves the minimax upper bound.

To estimate the probability matrix $M$, we propose using the Universal Singular Value Thresholding (USVT) method analyzed in \citet{chatterjee_matrix_2015} for estimating approximately low-rank matrices and later refined for networks in \citet{xu_rates_2018}.

\begin{algorithm}
\caption{USVT (\citet{chatterjee_matrix_2015,xu_rates_2018})}\label{alg:cap}
\begin{algorithmic}[1]
\State Input: Adjacency Matrix $A$, threshold $\tau > 0$
\State Compute the singular values of $A$, $s_1\geq s_2\geq \cdots \geq s_n$, where $\sum_{i=1}^n s_i u_i v_i^T$ is the singular value decomposition of $A$.
\State Collect those singular values above the threshold $S:=\{i:s_i \geq \tau\}$
\State Estimate $M$ by $\widehat{M}:= \sum_{i\in S}s_i u_i v_i^T,$ where we round to $0$ or $1$ if an entry happens to be less than zero or greater than 1 and zero out the diagonal.
\end{algorithmic}
\end{algorithm}

As discussed in \citet{chatterjee_matrix_2015}, the USVT method works well for estimating matrices that are low-rank or can be well-approximated by a low-rank matrix. We can give some intuition as to why graphons with spectral decay generate probability matrices that are approximately low rank. 

Suppose we have a graphon $W(x,y) = \sum_{i=1}^\infty \omega_i \varphi_i(x)\varphi_i(y)$ and an $n$-node network generated with latent node variables $\xi_1,...,\xi_n$. Set
$$
    \Phi_i = \begin{pmatrix}
        \varphi_i(\xi_1)\\
        \cdots\\
        \varphi_i(\xi_n)
    \end{pmatrix}.
$$

Then the probability matrix can be represented by $$M=\sum_{i=1}^\infty \omega_i \Phi_i \Phi_i^T - {\rm diag}\left(\sum_{i=1}^\infty \omega_i \Phi_i \Phi_i^T\right).$$

A natural low rank approximation is to use a truncated version 
$$
    M_k=\sum_{i=1}^k \omega_i \Phi_i \Phi_i^T,
$$
which has rank at most $k$. It can be seen here that as long as $\omega_i$ become small for $i > k$, the truncated probability matrix is a good low-rank approximation of the probability matrix $M$. This is the key intuition behind the effectiveness of USVT in this matrix estimation problem and is a vital part of the proof we present in the next section.

The USVT algorithm requires a choice of threshold $\tau$. We want to pick this threshold to be slightly above $\opnorm{A-M}$. Intuitively, this is because we want to discard the eigenvalues that are smaller than the largest eigenvalue of the noise since the since the signal in those eigencomponents are dominated by noise. Standard matrix concentration inequalities show that $\opnorm{A-M} \asymp \sqrt{n}$ and thus, we want to choose $\tau \asymp \sqrt{n}$. For implementation, the constant can be chosen to be $4$ based on the analysis given in \citet{xu_rates_2018}.

We now present the main result upper bounding the error of the USVT method.

\begin{theorem}
    \label{thm: upper bound}
 For any probability matrix $M$ generated by \eqref{eq:W2M} with a graphon $W\in\mathcal W(\alpha)$ and $\alpha > \frac{1}{2}$, the MSE for the USVT estimate has upper bound
    $$
    \frac{1}{n^2} E\left[\lfro{\widehat{M}-M}^2\right] \lesssim n^{-\frac{2\alpha - 1}{2\alpha}},
    $$
    where the expectation is taken over both $\xi$ and $A$.
\end{theorem}
Our proof of Theorem \ref{thm: upper bound} largely follows the scheme developed in \cite{xu_rates_2018} for upper bounding probability matrix estimation error for graphons with eigen decay. One technical difficulty in our case is controlling the discrepancy between the graphon and probability matrix incurred by the diagonal, which is substantially more straightforward when the graphon is smooth.  This part is detailed in the proof of Lemma \ref{lm: tail decay} below.

In addition to showing this upper bound, we also show that the rate of USVT is optimal up to log terms in this graphon class in the minimax sense. This is different in the case where we only consider SBM or H\"{o}lder smooth graphons \citet{gao_rate-optimal_2015,xu_rates_2018,klopp_oracle_2017}.

\begin{theorem}
    \label{thm: lower bound}
     Let $\xi_1,...,\xi_n \sim \text{Uniform}([0,1])$ and $M$ be the probability matrix induced by graphons in $\mathcal{W}(\alpha)$.  For $\alpha > 1$, we have
    $$
            \inf_{\widehat{M}} \sup_{\mathcal{W
            }(\alpha)} \frac{1}{n^2} E\left[\lfro{\widehat{M}-M}^2\right] \gtrsim (n\log(n))^{-\frac{2\alpha-1}{2\alpha}}\,.
    $$
\end{theorem}

\begin{remark}
    The upper and lower bounds have matching rates (up to a log factor) of convergence for $\alpha>1$.  The difficulty in proving matching lower bounds for $\alpha\le 1$ is in constructing a packing set of the intersection between the Stiefel manifold and the $\ell_\infty$ ball.  We conjecture the same minimax rate should still hold for a certain range of $\alpha$ smaller than 1 and consider this an open problem.
\end{remark}
\paragraph{Matching Rates} Our derivation of the minimax lower bound matches the upper bound up to a log factor. We discuss two settings in which we can get exactly matching rates.

In our graphon model of networks, the latent variables $\xi_1,...,\xi_n$ are sampled iid from the Uniform$([0,1])$ distribution. This is natural from the exchangeable random array and random graph perspective \citep{aldous_representations_1981,hoover_relations_1979,lovasz_large_2012}. Previous work in graphon estimation has studied graphon models where the latent variables can be sampled from arbitrary distributions. In this setting, we have matching upper and lower bounds.
\begin{theorem}
    \label{thm: lower bound arbitrary dist}
    Let $\mathcal{P}$ be the collection of all distributions on $[0,1]^n$. For $\alpha > 1$, we have
    $$
            \inf_{\widehat{M}} \sup_{\mathcal{W
            }(\alpha)} \sup_{P_\xi \in \mathcal{P}}\frac{1}{n^2} E\left[\lfro{\widehat{M}-M}^2\right] \asymp n^{-\frac{2\alpha-1}{2\alpha}}\,,
    $$
    where expectation is taken over $A$ and $(\xi_1,...,\xi_n) \sim P_\xi$.
\end{theorem}

The second setting where matching rates is possible relates to the packing of subspaces. Let $\mathbf{V}_{n,k}^\circ$ denote the set of orthonormal matrices $V \in \R^{n\times k}$ such that the columns all sum to zero.
\begin{assumption}
\label{as: pajor packing linf}
    There exists a collection $\mathcal{V} := \{V_1,...,V_N\} \subset \mathbf{V}_{n,k}^\circ$ satisfying
        \begin{enumerate}
        \item For $V,V' \in \mathcal{V}$, we have $\lfro{VV^T - V'{V'}^T}^2 \gtrsim k$
        \item $\log N \gtrsim nk$
        \item For any $V \in \mathcal{V},\linf{VV^T} \lesssim \frac{k^{\beta}}{n}$, where $0<\beta < \frac{1}{4}$.
    \end{enumerate}
\end{assumption}
Under this assumption, we can again achieve matching bounds.
\begin{theorem}
    \label{thm: lower bound linf packing}
    Suppose Assumption \ref{as: pajor packing linf} is satisfied. Then for $\alpha > 1$, we have
    $$
            \inf_{\widehat{M}} \sup_{\mathcal{W
            }(\alpha)} \frac{1}{n^2} E\left[\lfro{\widehat{M}-M}^2\right] \asymp n^{-\frac{2\alpha-1}{2\alpha}}\,.
    $$
\end{theorem}

The packing conditions given in Assumption \ref{as: pajor packing linf} is important in constructing the graphons used in the lower bound. The main difficulty is controlling the $\ell$-infinity norm $\linf{VV^T}$. In our proof, we are able to construct such a packing where $\linf{VV^T} \lesssim \frac{k}{n}$. This additional slack results in the additional log term. We conjecture that it may be possible to construct a packing set satisfying Assumption \ref{as: pajor packing linf} with more careful control of the $\ell$-infinity norm resulting in matching rates.

\section{Proofs}
\label{sec: proofs}
In this section, we present the proofs of our main results. To simplify the exposition we leave the proofs of the more technical lemmas to Appendix \ref{sec: main appendix}.

\subsection{Proof of Upper Bound (Theorem \ref{thm: upper bound})}

The performance of universal singular value thresholding for estimating probability matrices of networks has been thoroughly analyzed in \citet{xu_rates_2018}. For completeness, we summarize the argument in \citet{xu_rates_2018}, refining to account for our graphon spectral decay condition when necessary.

The main result in \citet{xu_rates_2018} is the following upper bound on the estimation error of USVT.

\begin{proposition}
\label{prop: USVT bound}(Theorem 1 \citep{xu_rates_2018})
    When the USVT threshold $\tau \asymp \sqrt{n}$, the estimation error of the USVT estimate is bounded above by
    \begin{equation*}
        \frac{1}{n^2}E\lfro{\widehat{M}-M}^2 \lesssim \min_{0\leq k\leq n}\left(\frac{k}{n} + \frac{1}{n^2}\sum_{i\geq k+1}E[\lambda_i^2]\right).
    \end{equation*}
\end{proposition}

The key takeaway is that the estimation error can be controlled by balancing $\frac{k}{n}$ the estimation error for estimating a rank $k$ matrix and $\frac{1}{n^2}\sum_{i\geq k+1}E[\lambda_i^2]$ the approximation error for approximating $M$ by a rank $k$ matrix.

The performance of this method depends on the decay of the tail eigenvalues of $M$. We can characterize this tail decay as follows.

\begin{definition}[Definition 1 \citep{xu_rates_2018}]
\label{def: matrix tail decay}
    We say the network probability matrix $M$ has polynomial tail decay with rate $\beta$ if
    $$
        \frac{1}{n^2} \sum_{i=k+1}^n E[\lambda_i^2]\leq c_0 k^{-\beta} + c_1 n^{-1},
    $$
    for all $1\leq k \leq n-1$ where $\lambda_i$ are the eigenvalues of $M$ in decreasing order in absolute value, the expectation is taken over the latent $\xi_i$ and $c_0,c_1$ are constants.
\end{definition}

When the probability matrix $M$ has polynomial tail decay with rate $\beta$, the approximation rate is approximately $k^{-\beta}$. Thus, the optimal trade-off between estimation and approximation rates occurs at $k \asymp n^{\frac{1}{\beta+1}}$, leading to the following upper bound.

\begin{proposition}[Corollary 1 \citep{xu_rates_2018}]
\label{prop: usvt rate}
    Let $M$ have polynomial tail decay with rate $\beta$ and $\tau \asymp \sqrt{n}$, then
    $$
        \frac{1}{n^2}E\lfro{\widehat{M}-M}^2 \lesssim n^{-\frac{\beta}{\beta + 1}}.
    $$
\end{proposition}

Thus, the key is to show how having decay condition \ref{def: spectral decay} on the eigenvalues of the graphon $W$ corresponds to tail decay of the corresponding probability matrix $M$ as in Definition \ref{def: matrix tail decay}. This correspondence is given in the following lemma.

\begin{lemma}
\label{lm: tail decay}
    If $M$ is the probability matrix of a graphon $W$ with polynomial spectral decay $\alpha$ with $\alpha > \frac{1}{2}$, then $M$ has tail polynomial spectral decay with rate $2\alpha-1$.
\end{lemma}

We save the proof of this result in the appendix, but we can give some simple intuition of the result. The term we are controlling corresponds to the approximation error of the probability matrix $M$ by a rank $k$ matrix. For intuition, let's work with graphons instead of probability matrices. The best rank $k$ approximation of a graphon $W(x,y) = \sum_{i=1}^\infty \omega_i \varphi_i(x)\varphi_i(y)$ is the graphon $W_k(x,y) = \sum_{i=1}^k \omega_i \varphi_i(x)\varphi_i(y).$ The error between this low-rank approximation is 
$$
    \ltwo{W-W_k}^2 = \sum_{i=k+1}^\infty \omega_i^2\,.
$$
When $W$ has polynomial spectral decay with rate $\alpha$, by using integral comparison, we can bound
$$
\sum_{i=k+1}^\infty \omega_i^2 \lesssim k^{-(2\alpha - 1)}\,,
$$
which shows the $2\alpha -1$ rate in the approximation error. The proof of Lemma \ref{lm: tail decay} will carry this intuition on the graphon level to the probability matrix level. The main difference is that at the probability matrix level, there are no self-edges, so there is an additional error corresponding to the self-edges (the diagonal of the graphon). Under the mild regularity provided by Assumption \ref{as: L_2 diagonal}, this error is of order $n^{-1}$ which satisfies the decay condition given in Definition \ref{def: matrix tail decay}.

The proof of Theorem \ref{thm: upper bound} is then an application of Lemma \ref{lm: tail decay} and Proposition \ref{prop: usvt rate}.

\begin{proof}(Theorem \ref{thm: upper bound})
    Let $W$ be a graphon with polynomial spectral decay $\alpha$. Then by Lemma \ref{lm: tail decay}, the probability matrix has polynomial tail decay with rate $2\alpha - 1$. Thus Proposition \ref{prop: usvt rate} shows that 
    \begin{align*}
    \frac{1}{n^2} E\left[\lfro{\widehat{M}-M}^2\right] \lesssim &n^{-\frac{2\alpha - 1}{2\alpha}},
    \end{align*}
    which shows the upper bound.
\end{proof}

\subsection{Proof of Lower Bound (Theorem \ref{thm: lower bound})}

Our proof of the lower bounds involves using the following formulation of Fano's lemma \citep{yu_assouad_1997} which was also used in previous work in graphon estimation \citep{gao_rate-optimal_2015,klopp_oracle_2017}.

Let $(\mathcal{X},\mathcal{A})$ be a measureable space of the observations and let $(\Theta, d)$ be a parameter space with pseudometric $d$. In this setting, each $\theta \in \Theta$ induces a probability measure $P_\theta$ on $(\mathcal{X},\mathcal{A})$. Let $T_N = \{\theta_1,...,\theta_N\} \subset \Theta$ be finite collection of parameters.

Fano's inequality gives the following lower bound

\begin{proposition}[Fano]
\label{prop: fano}
    Suppose that the parameters in $T_N$ are separated by $\alpha_N$, that is for $j\neq j'$
    $$
        d(\theta_j,\theta_{j'}) \geq \alpha_N
    $$
    and that 
    $$
        D_{KL}(P_{\theta_j},P_{\theta_j'}) \leq \beta_N.
    $$

    Then,
    $$
        \max_j E_{\theta_j} d(\hat{\theta},\theta_j) \geq \frac{\alpha_N}{2} \left(1-\frac{\beta_N + \log 2}{\log N}\right).
    $$
\end{proposition}

The probability estimation problem fits into the Fano framework in the following way. Let $I_i := \left[\frac{i-1}{m},\frac{i}{m}\right]$, for $i=1,...,m$, where $m := \frac{n}{4\log(n)}.$ For $i=1,...,m$, let $s_i$ denote the the number of $\xi_k \in I_i$, where $\xi_1,...,\xi_n \overset{iid}{\sim} \text{Uniform}([0,1]).$ To apply Fano's method, we first want to condition on a given realization of the latent $\xi$'s. The following lemma is useful for defining the conditioning set.
\begin{lemma}
    \label{lm: conditioning set}
   There exists positive constants $\lambda_1$, $\lambda_2$ such that 
   $$
    P\left(\frac{\lambda_1 n}{m} \leq s_i \leq \frac{\lambda_2 n}{m} \text{ for all } i=1,...,m\right) > c > 0,
   $$
   where $c$ is a constant independent of $n$.
\end{lemma}
 The idea is that when $\xi_j \overset{iid}{\sim}\text{Uniform}([0,1])$, each interval $I_i$ has around $\frac{n}{m}$ number of $\xi$'s with high probability when $n/m\gtrsim\log n$.

Let $E_n$ denote the event satisfying Lemma \ref{lm: conditioning set}. In our problem, the parameter space $\Theta$ consists of all probability matrices that can be generated by a graphon $W \in \mathcal{W}(\alpha)$ given the latent variables $\xi_1,...,\xi_n \in E_n$. That is, given a graphon $W$, let $M_W$ be the matrix such that $M_{W_{ij}} = W(\xi_i,\xi_j)$. Then $\Theta = \{M_W: W \in \mathcal{W}(\alpha)\}. $ For $M,M' \in \Theta$, we use the pseudometric
$$
    d(M,M') = \frac{1}{n^2} \lfro{M-M'}^2.
$$
Given $M\in \Theta$, let $P_{st} := \text{Bernoulli}(M_{st})$. Then $M$ parameterizes the distribution $\otimes_{s,t} P_{st}$.

The general idea of Fano's method is to reduce the problem to hypothesis testing over the discrete set $T_N$. The conditions given in Proposition \ref{prop: fano} give a general idea on how to choose the finite set $T_N$. We want a large collection (large $N$) which is well-separated in the pseudometric $d$ (large $\alpha_N$) but the induced probability measures are very similar in KL divergence (small $\beta_N$).

We now describe how to construct such a collection $T_N$. We do this by constructing graphons from subspaces of $\R^{m}$. Set $k\asymp (n\log(n))^{\frac{1}{2\alpha}}$.
Let $\mathbf{V}_{m,k}^\circ$ denote the set of orthonormal matrices in $\R^{m\times k}$ where the columns are orthogonal to $(1,...,1)^T$. The columns of $V \in \mathbf{V}_{m,k}^\circ$ define an orthonormal basis for a $k$-dimensional subspace of $\R^{m}$. 

Each $V \in \mathbf{V}_{m,k}^\circ$ forms an orthonormal set of functions in $L_2([0,1])$. Let $v_1,...,v_k$ denote the columns of $V$. For each $v_i, i=1,...,k$, define the function
$$
    \varphi_{v_i}(x) := \sqrt{m} \sum_{j=1}^{m} v_{ij} \mathds{1}_{(\frac{j-1}{m},\frac{j}{m}]}(x).
$$

A simple calculation shows that $\ltwo{\varphi_{v_i}}^2 = 1$ for $i=1,...,k$ and $\langle\varphi_{v_i},\varphi_{v_i'}\rangle_2 = 0$ for $i\neq i'$. Also, see that $\langle \varphi_{v_i}, \mathbb{1}\rangle = 0$, where $\mathds{1}$ is the constant function which has value $1$ everywhere. Thus, we can use these as eigenfunctions for a rank $k+1$ graphon. In particular, for $V \in \mathbf{V}^\circ_{m,k}$, we define the graphon
$$
W_V(x,y) = \frac{1}{2} + k^{-\alpha} L\sum_{i=1}^k \varphi_{v_i}(x)\varphi_{v_i}(y),
$$
where $L >0$ is a constant to be chosen later. This defines a class of graphons with polynomial spectral decay $\alpha$. We will denote the probability matrix induced by graphon $W_V$ by $M_V$.

To use Fano's method, we need to compute the KL-divergence and separation of probability matrices constructed from such graphons. Let's start with the KL-divergence. We use the following lemma given in \citet{klopp_oracle_2017}.

\begin{lemma}
    \label{lm: KL computation}
    Suppose that $P = \text{Bernoulli}(p)$ and $Q = \text{Bernoulli}(q)$. If $p,q \in \left[\frac{1}{4},\frac{3}{4}\right]$, then
    $$
        D_{KL}(P,Q) \leq \frac{16}{3}(q-p)^2.
    $$
\end{lemma}

For $M_V,M_{V'} \in \Theta$, we have
$$
    D_{KL}(P_{M_V},P_{M_{V'}}) \leq \frac{16}{3} \lfro{M_V-M_{V'}}^2,
$$
using Lemma \ref{lm: KL computation}. 

We see that both the pseudometric and KL-divergence depend on the term $\lfro{M_V-M_{V'}}^2$. When these matrices come from subspaces $V,V' \in \mathbf{V}_{m,k}^\circ$, we show this term actually related to a notion of subspace distance. 

\begin{lemma}
    \label{lm: graphon dist to subspace dist}
    Let $M_V,M_{V'}$ be probability matrices corresponding to subspaces $V,V' \in \mathbf{V}_{m,k}^\circ$. Then,
    $$
    \lfro{M_V-M_{V'}}^2 \asymp n^2 k^{-2\alpha} \lfro{VV^T - V'{V'}^T}^2.
    $$
\end{lemma}

Notice that $VV^T$ is the projection matrix to project to the subspace spanned by the columns of $V$. Similarly, $V'{V'}^T$ is the projection matrix to project to the subspace spanned by the columns of $V'$. The distance between the projection matrices $\lfro{VV^T - V'{V'}^T}^2$ is a natural measure of distance between the subspaces themselves. We use the following result on packings of subspaces with respect to the projection matrix metric. It can be derived using the result of \citet{pajor_metric_1998} along with Varshamov-Gilbert bound \citep{varshamov_rom_rubenovich_evaluation_1957,gilbert_comparison_1952}. These results are stated in more detail in Appendix \ref{sec: auxiliary appendix}.

\begin{lemma}
    \label{lm: subspace packing} There exists $\mathcal{V} = \{V_1,...,V_N\} \subset \mathbf{V}_{m,k}^\circ$ such that
    \begin{enumerate}
        \item For $V,V' \in \mathcal{V}$, we have $\lfro{VV^T - V'{V'}^T}^2 \gtrsim k$
        \item $\log N \gtrsim mk$
        \item For any $V \in \mathcal{V},\linf{VV^T} \lesssim \frac{k}{m}$.
    \end{enumerate}
\end{lemma}

Let $\mathcal{V}$ be constructed as in Lemma \ref{lm: subspace packing}. Then for Fano's method, we choose 
$$
    T = \{M_{V} : V \in \mathcal{V}\}.
$$
We first need to check that $W_V$ are well-defined graphons or equivalently, that the entries of $M_{W_V}$ are between $[0,1]$.

See that the entries of $M_V$ is bounded by
$$
    \frac{1}{2} + 4k^{-\alpha} L m \linf{VV^T} \lesssim \frac{1}{2} + 4Lk^{1-\alpha} = \frac{1}{2} + 4L(n\log(n))^{\frac{1}{2\alpha}- \frac{1}{2}}.
$$
For $\alpha > 1$, $(n\log(n))^{\frac{1}{2\alpha}- \frac{1}{2}}$ is smaller than constant order so by choosing $L$ small independent of $n$, we can make $\frac{1}{2} + 4k^{-\alpha} L m \linf{VV^T} < 1$.

To conclude our proof, we need to compute the values of $\alpha_N$, $\beta_N$, and $\log N$ which are used in Proposition \ref{prop: fano}.

To upper bound the KL-divergence, see that for any $V,V' \in \mathcal{V}$
\begin{equation*}
    \begin{split}
        KL(P_{M_V},P_{M_{V'}}) &\lesssim n^2 k^{-2\alpha} \lfro{VV^T - V'{V'}^T}^2\\
        &\lesssim  n^2 k^{-2\alpha} \left(\lfro{VV^T}^2 + \lfro{V'{V'}^T}^2\right)\\
        &\lesssim n^2k^{1-2\alpha}\\
        & \lesssim n^2(n\log(n))^{\frac{1}{2\alpha}-1}.
    \end{split}
\end{equation*}
The third inequality follows since the squared Frobenius norm of a projection matrix to a $k$-dimensional subspace is $k$.

Thus, for Fano's inequality, we can take $\beta_N \lesssim n^2(n\log(n))^{\frac{1}{2\alpha}-1}$.

To lower bound the pseudometric $d$, see that for $V,V' \in \mathcal{V}$, we have 
\begin{equation*}
    \begin{split}
        d(M_V,M_{V'}) &= \frac{1}{n^2} \lfro{M_V - M_{V'}}^2\\
        &\gtrsim k^{-2\alpha} \lfro{VV^T - V'{V'}^T}^2\\
        &\gtrsim k^{1-2\alpha}\\
        &= (n\log(n))^{-\frac{2\alpha-1}{2\alpha}}.
    \end{split}
\end{equation*}
Thus, in Fano's inequality, we can take $\alpha_N = (n\log(n))^{-\frac{2\alpha-1}{2\alpha}}.$

In addition, by Lemma \ref{lm: subspace packing}, $\log N \gtrsim mk = n^2(n\log(n))^{\frac{1}{2\alpha}-1}$. Thus, by choosing $L$ sufficiently small, we can make
$$
    \left(1-\frac{\beta_N + \log 2}{\log N}\right) > 0.5,
$$
which establishes the lower bound of $\alpha_N = (n\log(n))^{-\frac{2\alpha-1}{2\alpha}}$ up to a constant.

Applying Fano's inequality then gives
for given $\xi_1,...,\xi_n \in E_n$,
$$
    \frac{1}{n^2} E\left[\lfro{\widehat{M}-M}^2\right] \gtrsim (n\log(n))^{-\frac{2\alpha-1}{2\alpha}},
$$
where expectation is taken only over $A$. Then by taking expectation over $\xi$ and using Lemma \ref{lm: conditioning set}, we get the desired result.

\section{Discussion}
\label{sec: discussion}
In this work, we investigate the optimality of spectral methods in estimating the probability matrix of a network generated by a graphon. Previous work in this area investigated piecewise constant and H\"{o}lder smooth graphons. In these classes there is a computational-statistical gap where the computationally friendly universal singular value thresholding method is not minimax optimal, while the minimax optimal method is not computationally feasible. By viewing graphons as integral operators, we see that they have a natural spectral decomposition. We construct a class of graphons whose eigenvalues in this spectral decomposition decay fast. In this class, we find that universal singular value thresholding is minimax optimal up to a log term for estimating the probability matrix and give settings where we can derive exact minimax rates. 

Our results concern estimating the probability matrix induced by graphons. A possible future direction is to consider estimating the graphon itself in terms of a cut $L_2$ metric. That is, given a graphon $W$, we want an estimate of the entire graphon $\widehat{W}$ under the loss
$$
    \inf_h \int\int (\widehat{W}(h(x),h(y)) - W(x,y))^2 dxdy,
$$
where the infimum is taken over all measure preserving bijections. Such questions have been explored previously in the graphon estimation literature for stochastic block models and smooth graphons \citep{klopp_oracle_2017}. 

\begin{singlespace}
\medskip
\bibliographystyle{plainnat}
\bibliography{graphon.bib}
\end{singlespace}

\appendix

\section{Proof of Main Lemmas}
\label{sec: main appendix}
We give the proofs of the main lemma's used in the proofs of the upper and lower bounds.

\subsection{Upper Bound Lemma}

\begin{proof}[Proof of Lemma \ref{lm: tail decay}]
    This argument largely follows the proof of Theorem 5 in \citet{xu_rates_2018}, along with computing the sum of tail graphon eigenvalues and handling of the diagonal term.
    
    We can bound $\frac{1}{n^2} \sum_{i=k+1}^n E[\lambda_i^2(M)]$ using low-rank approximations of $M$. See that if $N$ is at most rank $k$, then 
$$
    \frac{1}{n^2} \sum_{i=k+1}^n E[\lambda_i^2(M)] \leq \frac{1}{n^2} E[\lfro{M-N}^2].
$$

This low-rank approximation $N$ can be constructed by thresholding the spectral decomposition of graphons. Using the spectral decomposition of graphons, see that
$$
M = \begin{bmatrix} 
    0 & \cdots & \sum_{i=1}^\infty \omega_i \varphi_i(\xi_1)\varphi_i(\xi_n)\\
    \vdots & \ddots & \\
   \sum_{i=1}^\infty \omega_i \varphi_i(\xi_n)\varphi_i(\xi_1) &        & 0 
    \end{bmatrix}.
$$

Thus, a natural way to get an approximation $N$ is to truncate $M$ at $k$ leaving only the first $k$ eigenvalues
$$
N = \begin{bmatrix} 
    \sum_{i=1}^k \omega_i \varphi_i(\xi_1)\varphi_i(\xi_1) & \cdots &\sum_{i=1}^k \omega_i \varphi_i(\xi_1)\varphi_i(\xi_n) \\
    \vdots & \ddots & \\
   \sum_{i=1}^k \omega_i \varphi_i(\xi_n)\varphi_i(\xi_1) &        & \sum_{i=1}^k \omega_i \varphi_i(\xi_n)\varphi_i(\xi_n) 
    \end{bmatrix}.
$$

First, see that $N$ is at most rank $k$. To see this, let 
$$
\Phi_i := \begin{pmatrix}
    \varphi_i(\xi_1)\\
    \cdots\\
    \varphi_i(\xi_n)
\end{pmatrix}.
$$

Then,
$$
N= \sum_{i=1}^k \omega_i \Phi_i \Phi_i^T,
$$
which is the sum of $k$ rank $1$ matrices. Then by subadditivity of matrix rank, it is at most rank $k.$

Now we can bound the expected Frobenius norm
$$ 
    \frac{1}{n^2} E[\lfro{M-N}^2].
$$
We can work component-wise. First, lets consider the entry $E[|M_{uv}-N_{uv}|^2]$, when $u \neq v.$

See that
\begin{equation*}
    \begin{split}
        E[|M_{uv}-N_{uv}|^2] &= E\left(\sum_{i=k+1}^\infty \omega_i \varphi_i(\xi_u)\varphi_i(\xi_v)\right)^2\\
        &= E\left(W(\xi_u,\xi_v) - \sum_{i=1}^k \omega_i \varphi_i(\xi_u)\varphi_i(\xi_v)\right)^2\\
        &= \ltwo{W(\xi_u,\xi_v) - \sum_{i=1}^k \omega_i \varphi_i(\xi_u)\varphi_i(\xi_v)}^2\\
        &\leq \left(\ltwo{W(\xi_u,\xi_v) - \sum_{i=1}^m \omega_i \varphi_i(\xi_u)\varphi_i(\xi_v)} + \ltwo{\sum_{i=k+1}^m \omega_i \varphi_i(\xi_u)\varphi_i(\xi_v)}\right)^2\\
        &\leq 2\ltwo{W(\xi_u,\xi_v) - \sum_{i=1}^m \omega_i \varphi_i(\xi_u)\varphi_i(\xi_v)}^2 + 2\ltwo{\sum_{i=k+1}^m \omega_i \varphi_i(\xi_u)\varphi_i(\xi_v)}^2,
    \end{split}
\end{equation*}

where the inequality holds for all $m > k$ and $\ltwo{\cdot}$ denotes the $L_2$ norm integrating over both $\xi_u$ and $\xi_v$.

Then we can take the limit to get
$$
    E[|M_{uv}-N_{uv}|^2] \leq 2\lim_m \ltwo{W(\xi_u,\xi_v) - \sum_{i=1}^m \omega_i \varphi_i(\xi_u)\varphi_i(\xi_v)}^2 + 2\lim_m \ltwo{\sum_{i=k+1}^m \omega_i \varphi_i(\xi_u)\varphi_i(\xi_v)}^2.
$$

We know that
$$
\lim_m \ltwo{W(\xi_u,\xi_v) - \sum_{i=1}^m \omega_i \varphi_i(\xi_u)\varphi_i(\xi_v)}^2 = 0,
$$

as 
$$
W(\xi_u,\xi_v) = \sum_{i=1}^\infty \omega_i \varphi_i(\xi_u)\varphi_i(\xi_v)
$$ where the convergence is in $L_2$.

The other term expands as
$$
\lim_m \ltwo{\sum_{i=k+1}^m \omega_i \varphi_i(\xi_u)\varphi_i(\xi_v)}^2 = \lim_m E\sum_{k+1}^m \omega_i^2\varphi_i^2(\xi_u)\varphi_i^2(\xi_v) + \lim_m E\sum_{i\neq j}\omega_i\omega_j \varphi_i(\xi_u)\varphi_i(\xi_v)\varphi_j(\xi_u)\varphi_j(\xi_v).
$$
We deal with the cross terms in the sum first.
We have
\begin{equation*}
    \begin{split}
        E\left(\sum_{i\neq j} \omega_i\omega_j \varphi_i(\xi_u)\varphi_i(\xi_v)\varphi_j(\xi_u)\varphi_j(\xi_v)\right) &= \sum_{i\neq j} \omega_i\omega_j E[\varphi_i(\xi_u)\varphi_j(\xi_u)]E[\varphi_i(\xi_v)\varphi_j(\xi_v)]\\
        &= 0,
    \end{split}
\end{equation*}
where the first equality comes from independence of $\xi_u$ and $\xi_v$ and the second equality comes from orthogonality of the eigenfunctions.

For the other term, we have that 
\begin{equation*}
    \begin{split}
        E\sum_{i=k+1}^m \omega_i^2 \varphi_i(\xi_u)^2\varphi_i(\xi_v)^2 &= \sum_{i=k+1}^m \omega_i^2 E[\varphi_i(\xi_u)^2]E[\varphi_i(\xi_v)^2] \\
        &= \sum_{i=k+1}^m \omega_i^2,
    \end{split}
\end{equation*}
where the first equality uses independence and the second uses that $L_2$ norm of eigenfunctions are 1.

Plugging everything back, we get
$$
E[|M_{uv}-N_{uv}|^2] \leq 2\sum_{i=k+1}^\infty \omega_i^2.
$$

Next we look at the diagonal terms, $M_{uu}$ for $u=1,...,n.$ See that
$$
E\left[|M_{uu}-N_{uu}|^2\right] = \norm{\sum_{i=1}^k \omega_i \varphi_i^2(\xi_u)}_2^2.
$$

This expands as
\begin{equation*}
    \begin{split}
         \norm{\sum_{i=1}^k \omega_i \varphi_i^2(\xi_u)}_2^2 & = \norm{W(\xi_u,\xi_u) - \sum_{i=k+1}^\infty \omega_i \varphi_i^2(\xi_u)}_2^2\\
         &\leq 2\norm{W(\xi_u,\xi_u)}_2^2 + 2\norm{\sum_{i=k+1}^\infty \omega_i \varphi_i^2(\xi_u)}_2^2.
    \end{split}
\end{equation*}

Let's work with these two terms separately. For the first term, see that
\begin{equation*}
    \begin{split}
        \norm{W(\xi_u,\xi_u)}_2^2 &= \norm{W(\xi_u,\xi_u)^2}_1 \\
        &\leq \norm{W(\xi_u,\xi_u)}_1 \norm{W(\xi_u,\xi_u)}_\infty.
    \end{split}
\end{equation*}

The inequality follows from H\"{o}lder. See that $\norm{W(\xi_u,\xi_u)}_\infty \leq 1$ because the range of $W$ is $[0,1]$. 

To handle the $\norm{W(\xi_u,\xi_u)}_1$ term, use that
\begin{equation*}
    \begin{split}
        \norm{W(\xi_u,\xi_u)}_1 &= \int \left|\sum_{i=1}^\infty \omega_i \varphi_i^2(\xi_u)\right| d\xi_u\\
        &\leq \int \sum_{i=1}^\infty |\omega_i|\varphi_i^2(\xi_u) d\xi_u\\
        &= \sum_{i=1}^\infty |\omega_i|\int \varphi_i^2(\xi_u)d\xi_u\\
        &= \sum_{i=1}^\infty |\omega_i| = \norm{W}_{\tr} < \infty.
    \end{split}
\end{equation*}

For the second term, under Assumption \ref{as: L_2 diagonal}, $\norm{\sum_{i=k+1}^\infty \omega_i \varphi_i^2(x)}_2$ viewed as a sequence in $k$ converges to $0$. In particular, there exists some $B$ independent of $k$ so that
$$
\norm{\sum_{i=k+1}^\infty \omega_i \varphi_i^2(\xi_u)}_2^2 < B,
$$
for all $k$. Note that $B$ is independent of $u$ as well since the norm doesn't depend on the specific $u$.

Then $E\left[|M_{uu}-N_{uu}|^2\right] \leq 2\norm{W}_{\tr} + 2B =: C$ which is a finite constant independent of $k$ and $u$.

Now using that there are $n$ diagonal terms and $n(n-1)$ off-diagonal terms, we can bound
\begin{equation*}
\begin{split}
\frac{1}{n^2} E[\lfro{M-N}^2]
 &\leq 2\sum_{i=k+1}^\infty \omega_i^2 + \frac{C}{n}.  
\end{split}
\end{equation*}

Now we can bound the first term by
    \begin{equation*}
        \begin{split}
            \sum_{i=k+1}^\infty \omega_i^2 & \leq \sum_{i=k+1}^\infty i^{-2\alpha}\\
            &\leq \int_{k}^\infty x^{-2\alpha} dx\\
            &= \frac{k^{-2\alpha+ 1}}{2\alpha - 1}.
        \end{split}
    \end{equation*}

Thus,
\begin{equation*}
\begin{split}
\frac{1}{n^2} \sum_{i=k+1}^n E[\lambda_i^2(M)] \leq \frac{1}{n^2} E[\lfro{M-N}^2]
 \leq \frac{2k^{-2\alpha + 1}}{2\alpha - 1} + \frac{C}{n},  
\end{split}
\end{equation*}

so $M$ has polynomial tail decay with rate $2\alpha-1$.
\end{proof}

\subsection{Lower Bound Lemmas}

In this section we prove the lemmas used in the lower bound proof.
\begin{proof}[Proof of Lemma \ref{lm: conditioning set}]
    We see that
    \begin{equation*}
        \begin{split}
        P\left(\bigcap_{i=1}^m \{s_i \geq \frac{\lambda n}{m}\}\right)
            &= 1 - P\left(\bigcup_{i=1}^m \{s_i < \frac{\lambda n}{m}\}\right)\\
            &\geq 1-\sum_{i=1}^m P\left(\{s_i < \frac{\lambda n}{m}\}\right).
        \end{split}
    \end{equation*}
    Let $X_k := \mathbb{1}(\xi_k \in I_1).$ We know that $X_k \sim \text{Bernoulli}(\frac{1}{m})$ and that $s_1 =\sum_{k=1}^n X_k.$ Now applying a Chernoff bound to sum of Bernoulli random variables such as \citet{chung_concentration_2006} Theorem 3.2, we see that
    \begin{equation*}
    \begin{split}
        P\left(s_1 < \frac{\lambda n}{m}\right) &\leq \exp\left(-\frac{n}{2m}(1-\lambda)^2)\right)\\
        &= \exp\left(-2\log(n)(1-\lambda)^2)\right)\\
        &= n^{-2(1-\lambda)^2},
    \end{split}
    \end{equation*}
    for $0\leq \lambda < 1$. Similarly, this holds for all $s_i$, $i=1,...,m$.
    Plug these back in to get
    $$
    P\left(\bigcap_{i=1}^m \{s_i \geq \frac{\lambda n}{m}\}\right) \geq 1- mn^{-(1-\lambda)^2} = 1-\frac{n^{1-2(1-\lambda)^2}}{4\log(n)}.
    $$ Choosing $\lambda = 1-\frac{1}{\sqrt{2}}$ gives that 
    $$
    P\left(\bigcap_{i=1}^m \{s_i \geq \frac{\lambda n}{m}\}\right) > 0.5
    $$
    for $n \geq 2$. A similar argument using the Chernoff bound for the upper tail shows that there exists constant $\lambda_2$ such that 
    $$
    P\left(\bigcap_{i=1}^m \{s_i \leq \frac{\lambda_2 n}{m}\}\right) > 0.5.
    $$
    Then as both events occur with probability larger than $0.5$, it must be that their intersection occurs with nonzero probability.
\end{proof}
\begin{proof}[Proof of Lemma \ref{lm: graphon dist to subspace dist}]
Define $\Phi_i = (\varphi_{v_i}(\xi_1),...,\varphi_{v_i}(\xi_n))^T$. Then
$$
    M_V = \frac{1}{2} + k^{-\alpha} L\sum_{i=1}^k \Phi_i \Phi_i^T.
$$
Similarly, we define $\Phi_i' = (\varphi_{v_i'}(\xi_1),...,\varphi_{v_i'}(\xi_n))^T$ and
$$
    M_{V'} = \frac{1}{2} + k^{-\alpha} L\sum_{i=1}^k \Phi_i' \Phi_i'^T.
$$
Then,
$$
    \lfro{M_V - M_{V'}}^2 = k^{-2\alpha} L^2 \lfro{\sum_{i=1}^k \Phi_i \Phi_i^T - \sum_{i=1}^k \Phi_i' \Phi_i'^T}^2.
$$
Let's show the lower bound. The upper bound follows from a similar argument. By construction of the conditioning set, we know that the entries of $\Phi_i$ contain $\sqrt{m}v_{ik}$ at least $\lambda_1 \frac{n}{m}$ times for $k=1,...,m$. Define
$$
    \Tilde{\Phi}_i := \sqrt{m}(\underbrace{v_{i1},...,v_{i1}}_{\lambda_1 \frac{n}{m}\text{times}},...,\underbrace{v_{im},...,v_{im}}_{\lambda_1 \frac{n}{m}\text{times}} )^T.
$$
Then see that
$$
    \lfro{\sum_{i=1}^k \Phi_i \Phi_i^T - \sum_{i=1}^k \Phi_i' \Phi_i'^T}^2 \geq \lfro{\sum_{i=1}^k \tilde{\Phi}_i \tilde{\Phi}_i^T - \sum_{i=1}^k \tilde{\Phi}_i' \tilde{\Phi}_i'^T}^2 = \lambda_1^2 n^2 \lfro{VV^T - V'V'^T}^2.
$$ The last equality follows as $\sum_{i=1}^k \tilde{\Phi}_i \tilde{\Phi}_i^T - \sum_{i=1}^k \tilde{\Phi}_i' \tilde{\Phi}_i'^T$ consists of $\lambda_1^2 \frac{n^2}{m^2}$ blocks of $VV^T-V'V'^T$. 
\end{proof}

\begin{proof}[Proof of Lemma \ref{lm: subspace packing}]
Let $m = m_1m_2,$ where $m_1 = 4k$. We first construct a packing of smaller matrices in $\mathbf{V}_{m_1,k}^\circ$. We will then combine these matrices into a packing of $\mathbf{V}_{m,k}^\circ$. The reason for using this two-stage procedure is to control the $\ell_\infty$ norm required in part 3 of the lemma.

For the first stage, we can apply Proposition \ref{prop: Pajor packing matrix}, to get a set $\{U_1,...,U_N\}$ of $m_1 \times k$ matrices satisfying
\begin{enumerate}
    \item $\lfro{UU^T-U'U'^T}^2 \gtrsim k$
    \item $\log N \gtrsim k(m_1-1-k) \gtrsim km_1.$
\end{enumerate}

We combine these smaller $m_1 \times k$ blocks to form the larger $m \times k$ matrices. We will use Varshamov-Gilbert bound to choose how to combine these matrices.

Let 
$$
    \Omega := \{\omega_1,...,\omega_M\} \subset \{1,...,N\}^{m_2},
$$

be a set of all N-ary sequences of length $m_2$, where 
$$
    \rho_H(\omega,\omega') \geq \frac{m_2}{4},
$$
where $\rho_H$ denotes Hamming distance (See Equation \ref{eq: Hamming Distance}). Varshamov-Gilbert (Proposition \ref{prop: vg bound}) gives a lower-bound on the size $M$
$$
    \log M \gtrsim m_2 \log N \gtrsim mk.
$$

For each $\omega \in \Omega$, define the $m \times k$ matrix
$$
    V_\omega := \begin{pmatrix}
        \frac{1}{\sqrt{m_2}}U_{\omega(1)}\\
        \cdots\\
        \frac{1}{\sqrt{m_2}}U_{\omega(m_2)}
    \end{pmatrix}.
$$

$V_\omega$ is an orthonormal matrix and by construction, the columns all sum to zero. Thus, $V_\omega \in \mathbf{V}_{m,k}^\circ$. We take $\mathcal{V} = \{V_\omega\}_{\omega \in \Omega}.$ The above analysis shows that $\log |\mathcal{V}| \gtrsim mk$ which shows part 2 of the lemma. Now we show part 1. That is,
    $$
        \lfro{V_\omega V_\omega^T - V_{\omega'} V_{\omega'}^T}^2 \gtrsim k,
    $$
    for all $\omega \neq \omega' \in \Omega.$

See that
    $$
    V_\omega V_\omega^T = \frac{1}{m_2}\begin{pmatrix}
        U_{\omega(1)}U_{\omega(1)}^T \cdots U_{\omega(1)}U_{\omega(n_2)}^T\\
        \cdots \cdots\\
        U_{\omega(n_2)}U_{\omega(1)}^T \cdots U_{\omega(n_2)}U_{\omega(n_2)}^T
    \end{pmatrix}.
$$
Then,
    \begin{equation*}
        \begin{split}
            \lfro{V_\omega V_\omega^T - V_\omega' V_\omega'^T}^2 &=\frac{1}{m_2^2} \sum_{ij}\lfro{U_{\omega(i)}U_{\omega(j)}^T-U_{\omega'(i)}U_{\omega'(j)}^T}^2.
        \end{split}
    \end{equation*}

We consider two cases for the individual terms in the sum. First, consider the block-diagonals. From the Pajor packing of $\mathbf{V}_{m_1,k}$, we know that
$$
    \lfro{U_{\omega(i)}U_{\omega(i)}^T-U_{\omega'(i)}U_{\omega'(i)}^T}^2 \gtrsim k,
$$
when $\omega(i) \neq \omega'(i)$ and is $0$ otherwise.

Now let's consider the off-diagonal blocks. By expanding, we have
\begin{align}
\label{eq: frob expansion}
        &\lfro{U_{\omega(i)}U_{\omega(j)}^T-U_{\omega'(i)}U_{\omega'(j)}^T}^2\nonumber \\
        =&  \lfro{U_{\omega(i)}U_{\omega(j)}^T}^2+\lfro{U_{\omega'(i)}U_{\omega'(j)}^T}^2 - 2tr\left(U_{\omega(i)}U_{\omega(j)}^TU_{\omega'(j)}U_{\omega'(i)}^T\right)\nonumber\\
        = & 2k - 2\tr\left(U_{\omega'(i)}^TU_{\omega(i)}U_{\omega(j)}^TU_{\omega'(j)}\right).
\end{align}
Next we upper bound the trace term $2\tr\left(U_{\omega'(i)}^TU_{\omega(i)}U_{\omega(j)}^TU_{\omega'(j)}\right)$. For notation, let $\theta_1,...,\theta_k$ denote the canonical angles between the subspaces defined by the columns of $U_{\omega(i)}$ and $U_{\omega'(i)}$ and $\tau_1,...,\tau_k$ be the canonical angles between the subspaces defined by the columns of $U_{\omega'(j)}$ and $U_{\omega(j)}$.

Then we know that $U_{\omega'(i)}^TU_{\omega(i)}$ has singular values $\cos(\theta_1),...,\cos(\theta_k)$ and $U_{\omega(j)}^TU_{\omega'(j)}$ has singular values $\cos(\tau_1),...,\cos(\tau_k)$.
Applying Von Neumann's trace inequality then gives
\begin{equation*}
\begin{split}
    |\tr\left(U_{\omega'(i)}^TU_{\omega(i)}U_{\omega(j)}^TU_{\omega'(j)}\right)| &\leq \sum_{\ell=1}^k \cos(\theta_\ell)\cos(\tau_\ell)\\ &\leq \sum_{\ell=1}^k \cos(\theta_\ell)^2 +\sum_{\ell=1}^k \cos(\tau_\ell)^2\\
    &= 2k - \sum_{\ell=1}^k \sin(\theta_\ell)^2 - \sum_{\ell=1}^k \sin(\tau_\ell)^2.
\end{split}
\end{equation*}

Now from Proposition \ref{prop: proj matrix distance}, we have that there exists constant $c>0$ such that
$$
 \sum_{\ell=1}^k \sin(\theta_\ell)^2 \geq c k\,~~\text{and~~}\sum_{\ell=1}^k \sin(\tau_\ell)^2 \geq c k,
$$ 
whenever $\omega(j)\neq \omega'(j)$ and $\omega(i)\neq \omega'(i)$,
so
$$
|\tr\left(U_{\omega'(i)}^TU_{\omega(i)}U_{\omega(j)}^TU_{\omega'(j)}\right)| \leq 2k - 2c k \,.
$$

Plug everything back into Equation \ref{eq: frob expansion} to get
$$
    \lfro{U_{\omega(i)}U_{\omega(j)}^T-U_{\omega'(i)}U_{\omega'(j)}^T}^2 \gtrsim k.
$$

Then,
    \begin{equation*}
        \begin{split}
            \lfro{V_\omega V_\omega^T - V_\omega' V_\omega'^T}^2 &=\frac{1}{m_2^2} \sum_{ij}\lfro{U_{\omega(i)}U_{\omega(j)}^T-U_{\omega'(i)}U_{\omega'(j)}^T}^2\\
            &\gtrsim \frac{1}{m_2^2} \rho_H^2(\omega,\omega') k\\
            &\gtrsim k.
        \end{split}
    \end{equation*}

Finally, we need to show part 3 of the Lemma.
For any $\omega$, by construction each row of $\sqrt{m_2}V_\omega$ is a row of orthonormal matrix, which as $\ell_2$-norm bounded by $1$. Thus each row of $V_\omega$ must have $\ell_2$-norm bounded by $1/\sqrt{m_2}$. Then the absolute value of the entries of $VV^T$ must be bounded by $1/m_2=4k/m$.

\end{proof}

\section{Proofs of Theorem \ref{thm: lower bound arbitrary dist} and \ref{thm: lower bound linf packing}}
In this section we give the proofs of the lower bound in the settings where matching rates can be achieved. The arguments are very similar to that of Theorem \ref{thm: lower bound}. Thus, for simplicity we only focus the major differences.

\begin{proof}[Proof of Theorem \ref{thm: lower bound arbitrary dist}] 

Let $P_\xi$ be the degenerate distribution such that $\xi_1=\frac{1}{n}, \xi_2=\frac{2}{n},...,\xi_n = \frac{n}{n}$. Set $k \asymp n^{\frac{1}{2\alpha}}$. By similar argument as in Lemma \ref{lm: subspace packing}, there exists $\mathcal{V} = \{V_1,...,V_N\} \subset \mathbf{V}_{n,k}^\circ$ such that
    \begin{enumerate}
        \item For $V,V' \in \mathcal{V}$, we have $\lfro{VV^T - V'{V'}^T}^2 \gtrsim k$
        \item $\log N \gtrsim nk$
        \item For any $V \in \mathcal{V},\linf{VV^T} \lesssim \frac{k}{n}$.
    \end{enumerate} For any $V \in \mathcal{V}$ define the graphon
    $$
    W_V(x,y) = \frac{1}{2} + k^{-\alpha} L\sum_{i=1}^k \varphi_{v_i}(x)\varphi_{v_i}(y),
    $$ in the same way as in the proof of Theorem \ref{thm: lower bound}. As $\xi_1=\frac{1}{n}, \xi_2=\frac{2}{n},...,\xi_n = \frac{n}{n}$, we can compute that
    $$
    \lfro{M_V-M_{V'}}^2 = n^2k^{-2\alpha} \lfro{VV^T-V'V'^T}^2 \asymp n^2k^{1-2\alpha},
    $$ for any $V,V' \in \mathcal{V}$. Then the same Fano computations as in the proof of Theorem \ref{thm: lower bound} gives the result.
\end{proof}

\begin{proof}[Proof of Theorem \ref{thm: lower bound linf packing}]
    Set $k \asymp n^{\frac{1}{2\alpha}}$ and let $\mathcal{V}$ be the packing given in Assumption \ref{as: pajor packing linf}. For any $V \in \mathcal{V}$ define the graphon
    $$
    W_V(x,y) = \frac{1}{2} + k^{-\alpha} L\sum_{i=1}^k \varphi_{v_i}(x)\varphi_{v_i}(y),
    $$ in the same way as in the proof of Theorem \ref{thm: lower bound}. Let $V,V' \in \mathcal{V}$ and denote their corresponding probability matrices by $M,M'$. Define 
    $$
    F(\xi_1,...,\xi_n) = \sum_{s,t} \left(\sum_{i=1}^k \varphi_{v_i}(\xi_s)\varphi_{v_i}(\xi_t)-\varphi_{v'_i}(\xi_s)\varphi_{v'_i}(\xi_t)\right)^2
    $$
    By construction, we have that
    $$
        \lfro{M-M'}^2 \asymp k^{-2\alpha}n^2F(\xi_1,...,\xi_n). 
    $$ See that $F$ satisfies bounded differences property with bound $\frac{2k^{2\beta}}{n}$. Applying McDiarmid's inequality, we get
    $$
        P\left(|F(\xi_1,...,\xi_n) - \lfro{VV^T-V'V'^T}^2| \geq \frac{c}{2}k\right) \leq 2\exp\left(-\frac{c^2k^{2(1-2\beta)}n}{8}\right),
    $$ where $c>0$ is a constant. Take union bound over all pairs $V,V'$, using there are $exp(2c'nk)$ such pairs by part two of Assumption \ref{as: pajor packing linf}, where $c'>0$ is a constant, we get
    \begin{align*}
        &P\left(|F(\xi_1,...,\xi_n) - \lfro{VV^T-V'V'^T}^2| \geq \frac{c}{2}k \text{ for all }V,V' \in \mathcal{V}\right)\\
        \leq & 2\exp\left(-\frac{c^2k^{2(1-2\beta)}n}{8}\right)\exp(2c'nk).
    \end{align*} As $1-2\beta >\frac{1}{2}$, the right side goes to $0$. Thus, there exists an event $E_n$, with $P(E_n) > C$ where $C>0$ is a constant independent of $n$ such that for all $\xi_1,...,\xi_n \in E_n$,
    $$
    |F(\xi_1,...,\xi_n) -\lfro{VV^T-V'V'^T}^2| \leq \frac{c}{2}k \text{ for all }V,V' \in \mathcal{V}.
    $$ In particular,
    $$
        \lfro{M-M'}^2 \asymp k^{1-2\alpha}n^2.
    $$ The remaining proof follows from similar Fano method computations.
\end{proof}
\label{sec: matching}
\section{Auxiliary Results}
\label{sec: auxiliary appendix}
The lower bound results use previous results regarding packing of subspaces and $m$-nary codes. For ease of reference we provide those results in this section and give references on more detailed exposition of these results.

\subsection{Subspace Packing}
Our construction of the lower bound uses techniques from subspace estimation. In this section, we give a brief outline of these tools following \citet{vu_minimax_2013}. A more detailed reference can be found in \citet{stewart_matrix_1990}.

We first define the metric that we will use to study the packing of subspaces. Let $\mathcal{E},\mathcal{F}$ be $k$-dimensional subspaces of $\R^n$. Let $E,F$ be the corresponding orthogonal projection matrices. We can define a notion of distance between subspaces by looking at their canonical angles.

\begin{definition}
    \label{def: cannonical angle}
    Let $s_1\geq \cdots \geq s_k$ be the non-zero singular values of $EF^\perp$. Then for $i=1,...,k$, define the $i$-th canonical angle between $\mathcal{E}$ and $\mathcal{F}$ to be
    $$
        \theta_i(\mathcal{E},\mathcal{F}) := \arcsin(s_i).
    $$
    We summarize all the canonical angles in a matrix
    $$
        \Theta(\mathcal{E},\mathcal{F}) = diag(\theta_1,...,\theta_k).
    $$
    Using this, we can define the distance between subspaces $\mathcal{E}$ and $\mathcal{F}$ to be
$
    \lfro{\sin(\mathcal{E},\mathcal{F}) \Theta}.
$
\end{definition}

Another natural candidate for a metric between subspaces is to consider a matrix norm between their orthogonal projection matrices. We will also use the following result which showing the relationship of this distance to the distance defined in Definition \ref{def: cannonical angle}.

\begin{proposition}[Proposition 2.1 \citep{vu_minimax_2013}]
\label{prop: proj matrix distance}
    Given subspaces $\mathcal{E},\mathcal{F}$ with orthogonal projection matrices $E,F$, we have
    $$
        \lfro{sin \Theta(\mathcal{E},\mathcal{F})}^2 = \frac{1}{2} \lfro{E-F}^2.
    $$
\end{proposition}

We are interested in constructing a large number of subspaces of $\R^n$ which are well-separated with respect to a subspace metric. Let $\mathbf{V}_{n,k}$ denote the set of all orthogonal $\R^{n\times k}$ matrices. In the literature, $\mathbf{V}_{n,k}$ is sometimes called the Stiefel manifold. Given $V \in \R^{n\times k}$, its columns define a $k$-dimensional subspace of $\R^n$. We will denote the corresponding subspace by $\mathcal{V}$. A vital part of the lower bound construction requires the studying the packing of subspaces $\mathcal{V}$ defined from $V \in \mathbf{V}_{n,k}$ in terms of the above defined subspace metric in Definition \ref{def: cannonical angle}. This packing was first studied by \citet{pajor_metric_1998}. The version we use is a case of Lemma A.6 in \citet{vu_minimax_2013}.
\begin{proposition}
\label{prop: Pajor packing}
    Suppose that $1 \leq k \leq n-k$. Then there exists a collection $\{V_1,...,V_N\} \subset \mathbf{V}_{n,k}$ such that
    \begin{enumerate}
        \item $\lfro{\sin \Theta(\mathcal{V}_i,\mathcal{V}_j)} \geq \sqrt{k}$ for $i \neq j$
        \item $\log N \gtrsim k(n-k)$.
    \end{enumerate}
\end{proposition}
Using Proposition \ref{prop: proj matrix distance}, we can rewrite part 1 as $\lfro{V_iV_i^T - V_jV_j^T}^2 \gtrsim k$ by noting that the orthogonal projection matrix to the column space of $V_i$ is given by $V_iV_i^T$.

We cannot directly apply Proposition \ref{prop: Pajor packing} for our purposes. The difficulty is that we want to study subspaces that are orthogonal to the space spanned by the vector with all entries $1$. To get around this issue, we can apply Pajor's packing on a smaller subspace. For notation, let $\R^{n\circ}$ denote the set of vectors $v \in \R^n$ such that $v \perp (1,...,1)^T$. Similarly, we use $\mathbf{V}_{n,k}^\circ$ to denote orthonormal $n\times k$ matrices where all the columns belong to $\R^{n\circ}$. Our goal is to construct a packing in $\mathbf V_{n,k}^\circ$. The space $\R^{n\circ}$ is isomorphic to $\R^{n-1}$. To make this explicit, by the Gram-Schmidt procedure, there exists an orthonormal basis of $\R^n, f_1,...,f_n$ such that $f_1,...,f_{n-1}$ is an orthonormal basis of $\R^{n\circ}$.

Now view $\R^{n\circ}$ as a copy of $\R^{n-1}$ with basis $f_1,...,f_{n-1}$. Applying the Pajor packing, there exists matrices $U_1,...,U_{N} \in \mathbf{V}_{n-1,k}$ such that $\lfro{U_iU_i^T - U_jU_j^T}^2 \gtrsim k$ and $\log N \gtrsim k(n-1-k)$. Let
$$
    V_i := \begin{bmatrix}
        f_1\,, & f_2\,, & \cdots & f_{n-1}
    \end{bmatrix} U_i.
$$
 See that
$$
    \lfro{U_iU_i^T - U_jU_j^T}^2 = \lfro{V_iV_i^T-V_j V_j^T}^2.
$$
See that $V_i \in \mathbf{V}_{n,k}^\circ.$ Thus, we have the following result
\begin{proposition}
\label{prop: Pajor packing matrix}
    Suppose that $1 \leq k \leq n-1-k$. Then there exists a collection $\{V_1,...,V_N\} \subset \mathbf{V}_{n,k}^\circ$, satisfying
    \begin{enumerate}
        \item $\lfro{V_iV_i^T - V_jV_j^T}^2 \gtrsim k$ for $i \neq j$
        \item $\log N \gtrsim k(n-1-k)$.
    \end{enumerate}
\end{proposition}

\subsection{N-nary Code Packing}
In statistical problems, the Varshamov-Gilbert bound for binary codes is often used. In our lower bound construction, we use a more general form of Varshamov-Gilbert for $m$-nary codes. For completeness, we summarize this result below. For reference, we follow the lecture notes \citep{guruswami_venkatesan_lecture_2010}.

For our purposes, we define an N-nary code of length $n$ to be a sequence taking values in $\{1,...,N\}$ of length $n$. For an N-nary code $\omega$, we use $\omega(i)$ to denote the $i$-th element of the sequence. We measure the distance between the sequences using Hamming distance. That is for two length $n$ N-nary codes $\omega,\omega'$, define the Hamming distance to be
\begin{equation}
\label{eq: Hamming Distance}
    \rho_H(\omega,\omega') := \sum_{i=1}^n \mathbb{1}(\omega(i) \neq \omega'(i)).
\end{equation}

Let 
$$
    \Omega := \{\omega_1,...,\omega_M\} \subset \{1,...,N\}^{n},
$$

be the set of all N-nary sequences of length $n$, where 
$$
    \rho_H(\omega_i,\omega_j) \geq \frac{n}{4},
$$
for $i \neq j$.
 We want a lower bound for $M$. Varshamov-Gilbert gives such a lower bound.

\begin{definition}
    \label{def: volume of Hamming ball}
    Denote by 
    $$
        V_N(n,r) := \sum_{j=0}^r \binom{n}{j} (N-1)^j
    $$
    the volume of Hamming ball of N-nary codes of length $n$ with radius $r$.
\end{definition}

Varshamov-Gilbert says that
\begin{lemma}[Varshamov-Gilbert]
    Let $M$ be the maximal size of N-nary codes of length $n$ with Hamming distance $d$ separation. Then we have the lower bound
    $$
        M \geq \frac{N^n}{V_N(n,d-1)}.
    $$
\end{lemma}

Applying Varshamov-Gilbert to our setting gives that
$$M \geq \frac{N^n}{V_N\left(n,\frac{n}{4}-1\right)} \geq \frac{N^n}{V_N\left(n,\frac{n}{4}\right).}$$.

We need to deal with the denominator. In this special case where the radius of the Hamming ball is a proportion of the length, this is easy to do. 

\begin{definition}
    Define the entropy function
    $$
        h_N(x) := x\log_N(N-1) -x\log_N x - (1-x)\log_N(1-x).
    $$
\end{definition}

Then we have the following bound
\begin{lemma}[Volume of Hamming ball]
\label{lm: hamming ball volume with entropy}
    For any $p \in [0,1-\frac{1}{N}]$,
    $$
        V_N(n,pn) \leq N^{h_N(p) n}.
    $$
\end{lemma}

We can apply this lemma and get the bound
$$
    M \geq N^{n\left(1-h_N\left(\frac{1}{4}\right)\right)},
$$
so
$$
    \log M \geq n\left(1-h_N\left(\frac{1}{4}\right)\right) \log N.
$$

One can check that for $N\geq 2$, we have 
$$
    0 < h_N\left(\frac{1}{4}\right) < 0.9.
$$

Thus, we get the following result.

\begin{proposition}
    \label{prop: vg bound}
    Let 
$$
    \Omega := \{\omega_1,...,\omega_M\} \subset \{1,...,N\}^{n},
$$

be the set of all N-nary sequences of length $n$, where 
$$
    \rho_H(\omega_i,\omega_j) \geq \frac{n}{4},
$$
for $i \neq j$. Then we have
$$
    \log M \gtrsim n \log N.
$$
\end{proposition}

\end{document}